\RequirePackage{ifpdf}
\ifpdf 
\documentclass[pdftex]{sigma}
\else
\documentclass{sigma}
\fi

\begin{document}


\renewcommand{\PaperNumber}{034}

\FirstPageHeading

\ShortArticleName{Natural and Projectively Invariant Quantizations on Supermanifolds}

\ArticleName{Natural and Projectively Invariant Quantizations\\ on Supermanifolds}

\Author{Thomas LEUTHER and Fabian RADOUX}

\AuthorNameForHeading{T.~Leuther and F.~Radoux}

\Address{Institute of Mathematics, Grande Traverse 12, B-4000  Li\`ege, Belgium}
\Email{\href{mailto:thomas.leuther@ulg.ac.be}{thomas.leuther@ulg.ac.be}, \href{mailto:fabian.radoux@ulg.ac.be}{fabian.radoux@ulg.ac.be}}

\ArticleDates{Received October 05, 2010, in f\/inal form March 23, 2011;  Published online March 31, 2011}

\Abstract{The existence of a natural and projectively invariant quantization in the sense of P.~Lecomte [{\em Progr. Theoret. Phys. Suppl.}   (2001), no.~144, 125--132] was proved by M.~Bordemann [math.DG/0208171], using the framework of Thomas--Whitehead connections. We extend the problem to the context of supermanifolds and adapt M.~Bordemann's method in order to solve it. The obtained quantization appears as the natural globalization of the $\mathfrak{pgl}({n+1|m})$-equivariant quantization on ${\mathbb{R}}^{n|m}$ constructed by P.~Mathonet and F.~Radoux in~[arXiv:1003.3320]. Our quantization is also a prolongation to arbitrary degree symbols of the projectively invariant quantization constructed by J.~George in~[arXiv:0909.5419] for symbols of degree two.}


\Classification{53B05; 53B10; 53D50; 58A50}

\section{Introduction}

The quantization is a concept that comes from physics. The quantization of a classical system whose phase space is a symplectic manifold $(\mathcal{M},\omega)$, consists in the construction of a Hilbert space~$H$ and a correspondence between classical and quantum observables. Classical observables are smooth functions on $\mathcal{M}$ while quantum observables are self-adjoint operators on~$H$. The concept of quantization has been formulated by P.~Dirac \cite{Dirac}, while trigged by the similarity between the formalisms of classical and quantum mechanics. At the beginning, the problem of the quantization consisted in f\/inding a linear bijection~$Q$ between $C^{\infty}(\mathcal{M})$ and a space of operators on~$H$ verifying three properties: the bijection~$Q$ has to transform the constant function~$1$ into the identity operator, the conjugation into the adjunction and the Poisson bracket into the commutator.

Prequantization \cite{Kirillov}, which associates with a function on $\mathcal{M}$ a dif\/ferential operator on the Hilbert space of complex functions square-integrable on $\mathcal{M}$, gives a positive answer to the Dirac problem. However, even when the phase space $\mathcal{M}$ is the cotangent bundle of a certain mani\-fold~$M$, it is not satisfactory because it gives us the Hilbert space $L^{2}(T^{*}M)$ which is too large for a~physically reasonable quantum system (it contains wave functions depending both on the position and the momentum coordinates, see~\cite[Chapter 5]{Woo} for details). Geometric quantization~\cite{Woo}, however, f\/ixes this issue by means of a polarization: if $\mathcal{M}=T^{*}M$, the vertical polarization allows of reducing the Hilbert space $L^{2}(T^{*}M)$ to the space~$L^{2}(M)$. Geometric quantization is then applied only to a subset of observables, those preserving the chosen pola\-ri\-za\-tion.

A priori, geometric quantization can be extended to the whole set of observables in many dif\/ferent ways. One can ask whether such an extension is unique if one requires extra conditions on the quantization map.

The uniqueness of a quantization procedure is often linked to a symmetry group. A vector f\/ield~$X$ on a manifold $M$ can be lifted in a natural way to a vector f\/ield on~$T^{*}M$, thereby def\/ining an action of the algebra of vector f\/ields on $M$, $\mathrm{Vect}(M)$, on the space of functions on $T^{*}{M}$ polynomial in the f\/ibers, called the space of symbols. It turns out that geometric quantization is the unique $\mathrm{Vect}(M)$-equivariant map from the space of symbols of degree less than or equal to one to the space of dif\/ferential operators on $M$, up to a normalization. However, geometric quantization cannot be extended to the whole space of symbols if one requires equivariance with respect to the Lie algebra $\mathrm{Vect}(M)$, due to cohomological reasons \cite{Cohomology}. One can ask whether there exists a Lie subalgebra $\mathfrak{g}\subset\mathrm{Vect}(M)$ for which the quantization map is $\mathfrak{g}$-equivariant. This~$\mathfrak{g}$ is supposed to be ``big enough'' to attain the uniqueness, but ``small enough'' to acquire the extension of the geometric quantization to the whole space of symbols.

When $M=\mathbb{R}^{n}$ with a $PGL(n+1,\mathbb{R})$-structure, the quantization map has been investigated by P.~Lecomte and V.~Ovsienko in \cite{LO}. They showed that there exists a unique quantization map that is $\mathfrak{pgl}(n+1,\mathbb{R})$-equivariant.

The concept of $\mathfrak{pgl}(n+1,{\mathbb{R}})$-equivariant quantization on $\mathbb{R}^{n}$ has a counterpart on an arbitrary manifold $M$ \cite{Leconj}. It aims at constructing a quantization map by means of a connection, depen\-ding only on its projective class (i.e.\ projectively invariant) and natural in all of its arguments. This quantization coincides with the projectively equivariant quantization when $M=\mathbb{R}^{n}$ and the connection is the projectively f\/lat one.

The existence of such a quantization procedure was proved by M.~Bordemann~\cite{Bor}. With
each class of torsion-free connections, he associated a unique linear connection, $\tilde{\nabla}$, on a principal line bundle $\tilde{M}\to M$. He managed to lift the space of symbols in a natural way to suitable tensors on $\tilde{M}$. He then applied the so-called \emph{standard ordering} to construct the projectively invariant quantization map.

Recently, P.~Mathonet and F.~Radoux \cite{MR5} constructed a super-version of the $\mathfrak{pgl}(n+1,{\mathbb{R}})$-equivariant quantization on ${\mathbb{R}}^{n}$. This super-quantization is a quantization on the super\-spa\-ce~${\mathbb{R}}^{n|m}$ which is equivariant with respect to a Lie superalgebra of supervector f\/ields on ${\mathbb{R}}^{n|m}$, this Lie superalgebra being isomorphic to $\mathfrak{pgl}(n+1|m)$.

In the same way as in the classical case, one can wonder if this super-quantization can be globalized. A partial positive answer to this question has been given in \cite{George}, where a
projectively invariant quantization on supermanifolds has been
constructed for symbols of degree less than or equal to two.

In this paper we prove the existence of a projectively invariant and natural quantization on a supermanifold. In Section~\ref{par:problem setting}, we recall the fundamental tools which we will need in the sequel (tensor densities, dif\/ferential operators, connections) and we give the def\/inition of a~natural projectively invariant quantization. Next, we solve the problem of existence of such a quantization by adapting Bordemann's method. To do this, we use the ``Thomas connection'' constructed by J.~George~\cite{George}, this connection being the super-version of the Thomas connection used by M.~Bordemann~\cite{Bor}.

\section{Problem setting}\label{par:problem setting}

For the sake of completeness, we brief\/ly recall in this section the def\/initions of tensor densities, dif\/ferential operators and symbols on supermanifolds. Then we set the problem of the existence of natural and projectively invariant quantizations on supermanifolds.
Throughout this paper, we denote by $M$ a smooth supermanifold of dimension $(n|m)$. The symbol $p(A)$ denotes the parity of the object $A$: it is equal to $0$ if $A$ is even, it is equal to $1$ if $A$ is odd.

Here we consider supermanifolds in the sense of F.~Berezin, B.~Kostant and D.~Leites (see \cite{Leites, MR1701618, MR1632008}).

\begin{definition}
A \emph{supermanifold} $M$ of dimension $(n|m)$ is a pair $(M_{0},\mathcal{A}_{M})$, where $M_{0}$ is a~$n$-dimensional smooth manifold and where $\mathcal{A}_{M}$ is a sheaf of superfunctions, i.e., a sheaf of superalgebras such that for all $p\in M_{0}$, there is an open neighbourhood $U\ni p$ together with an isomorphism
\[
\Phi_{U}: \ \mathcal{A}_{M}|_{U}\stackrel{\sim}{\longrightarrow} C^{\infty}_{{\mathbb{R}}^{n}}|_{U'}\otimes_{{\mathbb{R}}}\Lambda{\mathbb{R}}^{m},
\]
where $\Lambda{\mathbb{R}}^{m}$ denotes the exterior algebra of ${\mathbb{R}}^{m}$  and where $C^{\infty}_{{\mathbb{R}}^{n}}|_{U'}$ denotes the restriction to an open subset $U'$ of the sheaf of smooth functions on ${\mathbb{R}}^{n}$.

The images by $\Phi_{U}^{-1}$ of the canonical coordinates of ${\mathbb{R}}^{n}$ (resp. of the canonical generators of $\Lambda{\mathbb{R}}^{m}$) give even superfunctions $x^{1},\ldots,x^{n}$ (resp. odd superfunctions $x^{n+1},\ldots,x^{n+m}$). We call $(U,\Phi_{U})$ a local chart and $(x^{1},\ldots,x^{n},x^{n+1},\ldots,x^{n+m})$  \emph{graded coordinates} on~$M$.

 If $(U,\Psi_{U})$ is another local chart and if $(\bar{x}^{1},\ldots,\bar{x}^{n},\bar{x}^{n+1},\ldots,\bar{x}^{n+m})$ are corresponding graded coordinates, then the functions $\Phi_{U}(\bar{x}^{j})$, which form the \emph{change of graded coordinates}, will be denoted by $\bar{x}^{j}(x^{i})$.
\end{definition}
\begin{definition}
A \emph{morphism of supermanifolds} $M\to N$ is a pair $\phi=(\phi_{0},\phi^{*})$, where $\phi_{0}:M_{0}\to N_{0}$ is a smooth map and where $\phi^{*} : \mathcal{A}_{N} \to \mathcal{A}_{M}$ is a morphism of sheaves covering $\phi_{0}$, i.e., for every open subset $U\subset N_{0}$, $\phi^{*}$ gives a superalgebra morphism
\[\phi^{*}(U): \ \mathcal{A}_{N}(U)\to\mathcal{A}_{M}\big(\phi_{0}^{-1}(U)\big)\]
in a way compatible with the restriction maps. A \emph{diffeomorphism} is a morphism such that $\phi_{0}$ is a dif\/feomorphism and such that $\phi^{*}$ is an isomorphism.
\end{definition}

\subsection{Tensor densities}

The sheaf ${\mathcal F}_{\lambda,M}$ (or simply ${\mathcal F}_{\lambda}$) of $\lambda$-densities on $M$ is built from the Berezinian sheaf, whose formal def\/inition can be found in~\cite{Div}. For our purposes, it suf\/f\/ices to recall that, over an open set with graded coordinates $(x^{1},\ldots,x^{n},x^{n+1},\ldots,x^{n+m})$, a section of this sheaf, which we call a \emph{$\lambda$-density,} is expressed formally as
\[
\phi |Dx|^{\lambda},
\]
where $\phi$ is a local superfunction. Recall that, under a coordinate change ${\bar{x}}^i={\bar{x}}^i(x^j)$, the element $|Dx|^{\lambda}$ is multiplied by $|\mathrm{Ber}\, A|^{\lambda}$, where $\mathrm{Ber}$ denotes the Berezinian and where $A$ is the matrix corresponding to the change of coordinates, i.e., the matrix def\/ined by
\[
A_{j}^{i} = \frac {\partial{\bar{x}}^j}{\partial x^i}.
\]

\subsection{Dif\/ferential operators and symbols}

We denote by $\mathcal{D}_{\lambda,\mu,M}$ (or simply $\mathcal{D}_{\lambda,\mu}$) the sheaf of (f\/inite-order) \emph{differential operators} from $\lambda$-densities to $\mu$-densities.
For an open subset $U$ of $M_0$, elements $D \in \mathcal{D}_{\lambda,\mu}(U)$ are local ${\mathbb R}$-linear maps ${\mathcal F}_{\lambda}(U) \to {\mathcal F}_{\mu}(U)$ for which there is an integer $k \in {\mathbb N}$ such that in any system of local graded coordinates $(x^1,\dots,x^{n+m})$, the restriction of $D$ reads
\begin{gather}\label{eqDiffop}
\sum_{|\alpha|\leqslant k}D_\alpha \genfrac{(}{)}{}{}{\partial}{\partial x^1}^{\alpha_1}  \cdots  \genfrac{(}{)}{}{}{\partial}{\partial x^{n+m}}^{\alpha_{n+m}},
\end{gather}
where $\alpha$ is a multi-index, $|\alpha|=\sum\limits_{i=1}^{n+m}\alpha_i$, $\alpha_{n+1},\ldots,\alpha_{n+m}$ are in $\{0,1\}$ and $D_\alpha$ are local $\delta$-densities ($\delta = \mu - \lambda$).
The space $\mathcal{D}_{\lambda,\mu}(U)$ is thus f\/iltered by the order of dif\/ferential operators and we denote by $\mathcal{D}^k_{\lambda,\mu}(U)$ the space of dif\/ferential operators of order at most $k$.

As in the classical case, the highest order term of a dif\/ferential operator behaves, under a~change of graded coordinates, as a section of the sheaf of symbols, $\mathcal{S}^k_{\delta} = {\mathcal F}_{\delta}\otimes \vee^k {\mathcal T}_M$, where ${\mathcal T}_{M}$ denotes the tangent sheaf of~$M$ and where~$\vee$ denotes the supersymmetric product. This fact allows one to def\/ine the \emph{principal symbol operator} $\sigma_k:
\mathcal{D}^k_{\lambda,\mu}\to \mathcal{S}^k_{\delta}$. In graded coordinates,
\[
\sigma^k(D) =  \sum_{|\alpha|=k}D_\alpha (\partial_{1})^{\alpha_1}\vee\cdots\vee(\partial_{n+m})^{\alpha_{n+m}}
\]
if $D$ is given by~(\ref{eqDiffop}) and if $\partial_{1},\ldots,\partial_{n+m}$ denotes the canonical basis of local supervector f\/ields associated with the coordinate system $(x^{1},\ldots,x^{n+m})$. Moreover, we set $\mathcal{S}_{\delta}=\oplus_{k\geqslant 0}\mathcal{S}^{k}_{\delta}$.

\subsection{Superconnections and associated tensors}

A \emph{superconnection} on $M$ is a morphism of sheaves of super vector spaces $\nabla: {\mathcal T}_M \otimes_{\mathbb R} {\mathcal T}_M \to {\mathcal T}_M$ such that for any open subset $U$ in $M_{0}$, for any $X,Y \in {\mathcal T}_M(U)$ and any $f \in {\mathcal A}_M(U)$,
\[
\nabla_{fX} Y = f \nabla_X Y \qquad \text{while}  \quad \nabla_X {fY} = X(f) Y + (-1)^{p(X)p(f)} f \nabla_X Y,
\]
where $\nabla_X Y$ stands for $\nabla(X \otimes_{\mathbb R} Y)$.
Given graded coordinates $(x^1,\dots,x^{n+m})$ on $M$, a connec\-tion~$\nabla$ reads
\[
\nabla_X Y = \big( X^i \partial_i Y^k + (-1)^{p(i)(p(Y)+p(j))} X^i Y^j \Gamma_{ij}^k \big)\partial_k,
\]
where the so-called \emph{Christoffel symbols} $\Gamma_{ij}^k$ of $\nabla$ are superfunctions with parity $p(i)+p(j)+p(k)$.

From a superconnection $\nabla$ on $M$, one def\/ines its \emph{torsion} tensor ${\rm T}^{\nabla}$ by
\[
{\rm T}^{\nabla}(X,Y) = \nabla_X Y - (-1)^{p(X)p(Y)}\nabla_{Y}X - [X,Y].
\]
In graded coordinates, the vanishing of the torsion tensor of $\nabla$ translates into the supersymmetry of the Christof\/fel symbols in their low indices:
\[
\Gamma_{ij}^k = (-1)^{p(i)p(j)} \Gamma_{ji}^k.
\]
We denote by $\mathcal{C}_{M}$ (or simply $\mathcal{C}$) the sheaf of torsion-free superconnections on~$M$.

Remember also that from the curvature tensor ${\rm R}$ of $\nabla$, i.e.,
\[
{\rm R}(X,Y)Z = \nabla_X \nabla_Y Z - (-1)^{p(X)p(Y)} \nabla_Y \nabla_X Z - \nabla_{[X,Y]} Z ,
\]
the super-Ricci tensor ${\rm Ric}$ and the tensor ${\rm strR}$ are def\/ined as supertraces,
\begin{gather*}
{\rm Ric}(Z,Y)
  =
{\rm str}(X \mapsto(-1)^{p(X)(p(Y)+p(Z))}{\rm R}(X,Y)Z) ,
\\
{\rm strR}(X,Y)
  =
{\rm str}(Z \mapsto {\rm R}(X,Y)Z) .
\end{gather*}
The super-Ricci tensor ${\rm Ric}$ and the tensor ${\rm strR}$ are then given in coordinates by the following formulas:
\begin{gather*}
{\rm Ric}(Z,Y)
 =
(-1)^{p(i)(p(i)+p(Y)+p(Z))}({\rm R}(\partial_{i},Y)Z)^{i} ,
\\
{\rm strR}(X,Y)
 =
(-1)^{p(i)}({\rm R}(X,Y)\partial_{i})^{i},
\end{gather*}
if ${\rm R}(\partial_{i},Y)Z=({\rm R}(\partial_{i},Y)Z)^{j}\partial_{j}$ and ${\rm R}(X,Y)\partial_{i}=({\rm R}(X,Y)\partial_{i})^{j}\partial_{j}$.

\begin{remark}
Our tensor ${\rm strR}$, the supertrace of the curvature, appears for instance in \cite[Proposition 3.4]{Div}.
\end{remark}

\subsection{Projective equivalence of superconnections}\label{subsection: projective equivalence}

\begin{definition}
Two torsion-free superconnections $\nabla$, $\nabla'$ are \emph{projectively equivalent} if there is a~super $1$-form $\alpha$ on $M$ such that
\[
\nabla'_XY = \nabla_XY + \alpha(X) Y + (-1)^{p(X)p(Y)}\alpha(Y)X .
\]
An equivalence class for this relation is called a \emph{projective class}.
\end{definition}

Locally, the condition for two torsion-free superconnections $\nabla$ and $\nabla'$ to be projectively equivalent is expressed as $\Pi_{ij}^k = {\Pi'}_{ij}^k$, where
\begin{gather}\label{eq:fundamental descriptive invariant}
\Pi_{ij}^k = \Gamma_{ij}^k - \frac{1}{n-m+1}\big( \Gamma_{is}^s \delta_j^k (-1)^{p(s)} + \Gamma_{js}^s \delta_i^k (-1)^{p(i)p(j)+p(s)} \big) .
\end{gather}
Obviously, this characterization fails when the \emph{superdimension} $n-m$ is equal to $-1$ (in this case, the quantities $\Pi_{ij}^k$ cannot be def\/ined). We believe that a detailed study of projective superconnections could provide a geometric meaning for this singular situation, but this is yet to be done.

\subsection{Problem setting}

A \emph{quantization} on $M$ is an isomorphism $Q_M$ of sheaves of vector spaces, from the sheaf of symbols $\mathcal{S}_{\delta,M}$ to the sheaf of dif\/ferential operators $\mathcal{D}_{\lambda,\mu,M}$ such that for any $k \in {\mathbb
N}$ and any section $S$ of ${\mathcal S}^k_{\delta, M}$,
\[
\sigma^{k}(Q_M(S)) = S.
\]
Roughly speaking, a \emph{natural quantization} is a quantization which depends on a torsion-free superconnection and commutes with the action of superdif\/feomorphisms. More precisely, a~natural quantization is a collection of morphisms (def\/ined for every supermanifold $M$)
\[
Q_M:  \ \mathcal{C}_{M} \times \mathcal{S}_{\delta,M}\to \mathcal{D}_{\lambda,\mu,M}
\]
such that
\begin{itemize}\itemsep=0pt
\item for any section $\nabla$ of $\mathcal{C}_{M}$, $Q_M(\nabla)$ is a quantization;
\item if $\phi: M \to N$ is a local dif\/feomorphism of supermanifolds, then, for any section  $\nabla$ of $\mathcal{C}_{N}$ and any section $S$ of $\mathcal{S}_{\delta,N}$,
\[
Q_M(\phi^*\nabla)(\phi^*S) = \phi^*(Q_N(\nabla)(S)).
\]
\end{itemize}
A quantization $Q_M$ is \emph{projectively invariant} if one has $Q_M(\nabla) = Q_M(\nabla')$ whenever $\nabla$ and $\nabla'$ are projectively equivalent torsion-free superconnections.

The existence of natural and projectively invariant quantization is related to the existence of a $\mathfrak{pgl}(n+1|m)$-equivariant quantization in the sense of \cite{MR5} in the f\/lat situation.
\begin{definition}
When $n-m \neq -1$, we def\/ine the numbers
\[\gamma_{2k-l} = \frac{(n-m+2 k - l -(n-m+1)\delta)}{n-m+1}.\]
In this case, a value of $\delta$ is said to be \emph{critical} if there is $k,l\in {\mathbb{N}}$ such that $1\leq l\leq k$ and $\gamma_{2k-l}=0$. Notice that, in opposition to the classical situation, the value $\delta = 0$ can be critical since $n-m$ can be negative.
\end{definition}
One of the results of \cite{MR5} is then the following.
\begin{theorem}\label{plat1}
If $n-m \neq -1$ and $\delta$ is not critical, there exists a unique $\mathfrak{pgl}(n+1| m)$-equivariant
quantization. If $n-m = -1$, there exists a one-parameter family of $\mathfrak{pgl}(n+1| m)$-equivariant quantizations.
\end{theorem}

As in the classical context, the projective class of the canonical f\/lat superconnection is preserved by the vector f\/ields of $\mathfrak{pgl}(n+1|m)$. Therefore, if we can construct a natural projectively invariant quantization for supermanifolds of dimension $(n|m)$, it will necessarily coincide with the projectively equivariant quantization constructed in~\cite{MR5} when $M=\mathbb{R}^{n|m}$ and $\nabla$ is the projectively f\/lat superconnection.

\section{Thomas f\/iber bundle and connection}

As in the classical case, a projective class of torsion-free superconnections on a supermanifold def\/ines a superconnection on the associated ``Thomas bundle''. This fact was pointed out in \cite{George} and is brief\/ly recalled here.

\subsection{Thomas f\/iber bundle}\label{subsection: Thomas fiber bundle}

From a supermanifold $M$ of dimension $(n| m)$, the associated  supermanifold $\tilde{M}$ is constructed by adding an even coordinate $x^0$ to each coordinate system $(x^1,\dots ,x^{n+m})$. Under a change of coordinates ${\bar{x}}^i = {\bar{x}}^i(x^j)$ on $M$, this extra coordinate transforms according to the rule ${\bar{x}}^0 = x^0 + \log|{\rm Ber}\, A|$ where $A$ is the matrix corresponding to the change of coordinates.

By analogy with the classical situation, we set ${\mathcal E} = \partial_0$ and call it the \emph{Euler supervector field} of~$\tilde{M}$. The fact that ${\mathcal E}$ is well-def\/ined is easily seen from the transformation law of the components of a supervector f\/ield $X$ under a change of coordinates ${\bar{x}}^i = {\bar{x}}^i(x^j)$, namely
\[
X^{i}\partial_{x^i}=X^{i}\frac {\partial{\bar{x}}^j}{\partial x^i}\partial_{{\bar{x}}^j}.
\]

Densities on $M$ identify with some superfunctions on $\tilde{M}$. More precisely, we can associate with a $\lambda$-density expressed locally as $f = \phi |Dx|^\lambda$ the superfunction $\tilde{f}$ given by
\begin{gather}\label{formula: densities to superfunctions}
\tilde{f}\big(x^{0},x^{1},\ldots,x^{n+m}\big) = \phi\big(x^{1},\ldots,x^{n+m}\big)\exp\big(\lambda x^{0}\big).
\end{gather}
It follows directly from the transformation law of densities that $\tilde{f}$ is well-def\/ined. Moreover, it is $\lambda$-equivariant in the sense that
\[
{\rm L}_{\mathcal E} \tilde{f} = \lambda \tilde{f}.
\]
Conversely, from a $\lambda$-equivariant superfunction $\varphi$ on $\tilde{M}$, one def\/ines a $\lambda$-density $\varphi_0 |Dx|^\lambda$ on $M$ by setting
\begin{gather}\label{formula: superfunctions to densities}
\varphi_{0}\big(x^1,\dots ,x^n\big) = \varphi\big(x^0,x^1,\dots ,x^n\big)\exp\big(-\lambda x^0\big)
\end{gather}
for an arbitrary $x^0$. Because of the equivariance property of $\varphi$, the derivative of $\varphi_0$ with respect to $x^0$ is zero and the density is well-def\/ined. This way, we establish a one-to-one correspondence between $\lambda$-densities on~$M$ and $\lambda$-equivariant superfunctions on~$\tilde{M}$.

\subsection{Thomas connection}\label{subsection:Thomas connection}


The quantities $\Pi_{ij}^k$ given by (\ref{eq:fundamental descriptive invariant}) def\/ine the so-called \emph{fundamental descriptive invariant} of the projective class of $\nabla$. They can be used to construct the projectively invariant lift $\tilde{\nabla}$ of $\nabla$ to the supermanifold $\tilde{M}$. More precisely, $\tilde{\nabla}$ is def\/ined as the torsion-free superconnection on $\tilde{M}$ with Christof\/fel symbols
\begin{gather*}
\tilde{\Gamma}_{ij}^k = \Pi_{ij}^k , \qquad \tilde{\Gamma}_{0{\mathfrak a}}^{\mathfrak c} = \tilde{\Gamma}_{{\mathfrak a}0}^{\mathfrak c} =  \frac{-\delta_{\mathfrak a}^{\mathfrak c}}{n-m+1},
\\
\tilde{\Gamma}_{ij}^0 = \frac{n-m+1}{n-m-1} \big( \partial_q \Pi_{ij}^q - \Pi_{qi}^p \Pi_{pj}^q \big) (-1)^{p(q)(p(q)+p(i)+p(j))},
\end{gather*}
where $i$, $j$, $k$ ranges from $1$ to $n+m$ while ${\mathfrak a}$, ${\mathfrak c}$ ranges from~$0$ to~$n+m$.
Besides the singular case $n-m = -1$ (in which the quantities $\Pi_{ij}^k$ themselves cannot be def\/ined), one must also assume here that $n-m \neq 1$ in order to give sense to the above formulas. The latter hypothesis is the super analogue of the fact that M.~Bordemann's construction fails for a $1$-dimensional smooth manifold.

We shall now give a useful coordinate-free description of the lifted superconnection $\tilde{\nabla}$ in terms of the Euler supervector f\/ield of $\tilde{M}$ and horizontal lifts of supervector f\/ields on $M$. Using our coordinate system $(x^0,x^1,\dots ,x^{n+m})$ (cf.~Section~\ref{subsection: Thomas fiber bundle}), we introduce the so-called \emph{horizontal lift} to~$\tilde{M}$ of a supervector f\/ield $X = X^i \partial_i$ on $M$ by setting
\begin{gather}\label{eq: horizontal lift of vector fields}
 X^h = - (-1)^{p(s)} X^i \Gamma_{is}^s \partial_0 + X^i \partial_i .
\end{gather}
The vector f\/ield $X^{h}$ is well-def\/ined. Indeed, the derivative of the Berezinian of the matrix $A$ representing the change of coordinates can be computed in the following way \cite{Leites}:
\[
\partial_{x^{i}}(\mathrm{Ber}\;A)=(\mathrm{Ber}\;A)\mathrm{str}\big((\partial_{x^{i}}A)A^{-1}\big) ,
\]
where $(\partial_{x^i}A)_{l}^{k} = (-1)^{p(i)p(k)}\partial_{x^i}A_{l}^{k}$ and where
\[\mathrm{str}(B) = \sum_{i=1}^{p+q}(-1)^{p(i)(p(B)+p(i))}B_{i}^{i}\]
if $B\in\mathfrak{gl}(n|m)$.
We can then easily conclude using the transformation law of the Christof\/fel symbols:
\[
\bar{\Gamma}_{ij}^{k}=(-1)^{p(t)(p(l)+p(j))}
\left(-\frac {\partial x^{t}}{\partial {\bar{x}}^{i}}\frac {\partial x^{l}}{\partial {\bar{x}}^{j}}\frac {\partial^{2} {\bar{x}}^{k}}{\partial x^{t}\partial x^{l}}+\frac {\partial x^t}{\partial {\bar{x}}^i}\frac {\partial x^l}{\partial {\bar{x}}^j}\Gamma_{tl}^{r}\frac {\partial {\bar{x}}^k}{\partial x^{r}}\right).
\]

\begin{definition}\label{r}
We denote by $r$ the following multiple of a supersymmetric part of the Ricci tensor of $\nabla$:
\[
{\rm r}(X,Y) = \frac {1}{2(n-m-1)}\big({\rm Ric}(Y,X)+(-1)^{p(X)p(Y)}{\rm Ric}(X,Y)\big).
\]
\end{definition}

Finally, we are in position to give the coordinate-free description of the lifted supercon\-nec\-tion~$\tilde{\nabla}$. This description is useful in order to simplify computations in the sequel.

{\samepage
\begin{proposition}
With notations of \eqref{formula: densities to superfunctions}, we have for any supervector fields $X$, $Y$ on $M$,
\begin{gather*}
\tilde{\nabla}_{X^h}{Y^h}
= \left( \nabla_X Y \right)^h - \frac{1}{2} \widetilde{{\rm strR}(X,Y)}{\mathcal E} + (n-m+1)\widetilde{{\rm r}(X,Y)} {\mathcal E} ,
\\
\tilde{\nabla}_{X^h}{{\mathcal E}} = \tilde{\nabla}_{{\mathcal E}}{X^h} = \frac{-1}{n-m+1} X^h
, \qquad
\tilde{\nabla}_{{\mathcal E}}{{\mathcal E}} = \frac{-1}{n-m+1}{\mathcal E} .
\end{gather*}
\end{proposition}

\begin{proof}
The result is obtained after long but straightforward computations.
\end{proof}}

The lifted connection $\tilde{\nabla}$ is associated in a natural way with the connection $\nabla$ on $M$. Moreover, $\tilde{\nabla}$ is such that ${\rm L}_{\mathcal{E}}\tilde{\nabla}=0$, where
\[
{\rm L}_{\mathcal{E}}\tilde{\nabla}(X,Y) = [\mathcal{E},\tilde{\nabla}_{X}Y]-\tilde{\nabla}_{[\mathcal{E},X]}Y-\tilde{\nabla}_{X}[\mathcal{E},Y] .
\]
This invariance is due to the invariance of $\mathcal{E}$, of the horizontal lifts and of the functions $\widetilde{{\rm strR}(X,Y)}$ and $\widetilde{{\rm r}(X,Y)}$.

\section{Lift of symbols}
Our goal in this section is to lift in a natural and projectively invariant way a symbol~$S$ on~$M$ to a tensor $\tilde{S}$ on $\tilde{M}$. To this aim, we def\/ine in a f\/irst step a horizontal lift of~$S$ via the horizontal lift of supervector f\/ields~(\ref{eq: horizontal lift of vector fields}). In a second step, we def\/ine a map which transforms equivariant tensors on~$\tilde{M}$ into symbols on~$M$. We prove that the restriction of this map to the divergence-free tensors (with respect to~$\tilde{\nabla}$) is a bijection. The natural and projectively invariant lift is then the inverse map of this ``descent'' application.

\subsection{Horizontal lift of symbols}

Since a symbol $S$ of degree $k$ on $M$ is locally a sum of terms of the form $\phi| Dx|^{\delta} \otimes \partial_{i_1}\vee\cdots\vee\partial_{i_{k}}$, it suf\/f\/ices to def\/ine the horizontal lift on symbols of this form and to extend it by linearity.
\begin{definition}
The horizontal lift of a symbol $\phi| Dx|^{\delta}\otimes\partial_{i_1}\vee\cdots\vee\partial_{i_{k}}$, denoted by $(\phi| Dx|^{\delta}\otimes\partial_{i_1}\vee\cdots\vee\partial_{i_{k}})^{h}$, is given in our coordinate system on $\tilde{M}$ by the tensor $\tilde{f} \otimes \partial_{i_1}^{h}\vee\cdots\vee\partial_{i_{k}}^{h}$ if $f = \phi| Dx|^{\delta}$.
\end{definition}

{\sloppy We can easily observe that the horizontal lift of a symbol $S$ is $\delta$-equivariant, i.e.\ that \mbox{${\rm L}_{\mathcal{E}}S^{h}=\delta S^{h}$}. In particular, the horizontal lift of a $\delta$-density on $M$ to a superfunction on~$\tilde{M}$ coincides with the correspondence given in~(\ref{formula: densities to superfunctions}).

}

\subsection{Descent map}

Using the fact that a tensor of degree $k$ on $\tilde{M}$ can be locally decomposed in the basis $\partial_{1}^{h},{\ldots},\partial_{n+m}^{h}$, ${\mathcal E}$, it is enough to def\/ine the descent map on a tensor of the form
\begin{gather}\label{dec}
S=\sum_{l=0}^{k}\sum_{i_{1},\ldots,i_{k-l}}\varphi^{i_{1},\ldots,i_{k-l},0,\ldots,0}\otimes\partial_{i_1}^{h}\vee\cdots\vee\partial_{i_{k-l}}^{h}\vee{\mathcal E}^{l}.
\end{gather}
In the sequel, we denote by $\mathcal{S}^{k,\delta}_{\tilde{M}}$ the sheaf of $\delta$-equivariant tensors of degree $k$ on $\tilde{M}$.

\begin{proposition}
The map $\Psi$ whose value on a section $S$ of $\mathcal{S}^{k,\delta}_{\tilde{M}}$ expressed as in \eqref{dec}
 is given by
\[
\Psi(S)=\sum_{i_{1},\ldots,i_{k}} \varphi_0^{i_{1},\ldots,i_{k}}|Dx|^{\delta}\otimes \partial_{i_1}\vee\cdots\vee\partial_{i_{k}} ,
\]
where the coefficient $\varphi_0^{i_{1},\ldots,i_{k}}$ is given by \eqref{formula: superfunctions to densities}, is well-defined.
\end{proposition}

\begin{proof}
We can easily see that the form of $\Psi$ is preserved under a change of coordinates ${\bar{x}}^i = {\bar{x}}^i(x^j)$ on $M$. It is among other things due to the fact that
\begin{gather*}
\varphi^{i_{1},\ldots,i_{k}}\big({\bar{x}}^0,{\bar{x}}^1,\ldots,{\bar{x}}^{n}\big)\exp\big(-\delta {\bar{x}}^0\big)
  =  \varphi_0^{i_{1},\ldots,i_{k}}\big(x^{1},\ldots,x^{n}\big)|\mathrm{Ber}\,A|^{-\delta}
  =   \bar{\varphi}_0^{i_{1},\ldots,i_{k}}\big({\bar{x}}^1,\ldots,{\bar{x}}^{n}\big) ,
\end{gather*}
where $A$ stands for the matrix of the change of graded coordinates on $M$.
\end{proof}

It is possible to show that the map $\Psi$ is surjective exactly in the same way as in~\cite{Bor}. If $A_{k}$ is a symbol of degree~$k$ on~$M$, then the
sections of~$\mathcal{S}^{k,\delta}_{\tilde{M}}$ whose images by $\Psi$ are equal to~$A_{k}$ are those of the form
\begin{gather}\label{decomp}
A_{k}^{h}+A_{k-1}^{h}\vee{\mathcal E}+\cdots+A_{0}^{h}\vee{\mathcal E}^{k},
\end{gather}
for some sections $A_{k-j}$ of $\mathcal{S}_{\delta}^{k-j}$ for $j=1,\ldots,k$.

\subsection{Projectively invariant lift of symbols}

In the sequel, we denote by $dx^{1},\ldots,dx^{n+m}$ the dual basis of the canonical basis of local supervector f\/ields $\partial_{1},\ldots,\partial_{n+m}$ on $M$. It means that $dx^{i}(\partial_{j})=\delta_{j}^{i}$ for all $i$, $j$.

The covariant derivative with respect to $\nabla$ of a $\delta$-density $\phi| Dx| ^{\delta}$ in the direction of a supervector f\/ield $X$ is given in coordinates by
\[
\nabla_{X}\big(\phi| Dx| ^{\delta}\big) = \big( X.\phi-(-1)^{p(s)}\delta X^{i}\Gamma_{is}^{s}\phi \big) | Dx| ^{\delta}.
\]
We can easily show that this formula is preserved under a change of coordinates and that it def\/ines a covariant derivative such that
\[
\widetilde{\nabla_{\partial_{i}}\phi| Dx| ^{\delta}} = \partial_{i}^{h} . \widetilde{\big( \phi| Dx| ^{\delta}\big)} .
\]
\begin{definition}\label{interior}
The interior product of a symbol $S = \phi| Dx|^{\delta} \otimes \partial_{i_1}\vee\cdots\vee\partial_{i_{k}}$ by a super $1$-form $\alpha=\alpha_{i} dx^{i}$ is def\/ined by
\[
{\rm i}(\alpha)(S) =\sum_{j=1}^{k}(-1)^{p(\alpha)(p(\phi)+p(i_{1})+\cdots+p(i_{j-1}))}\phi| Dx|^{\delta}\otimes\partial_{i_1}\vee\cdots\vee\overset{(j)}{\alpha_{i_{j}}}\vee \cdots\vee\partial_{i_{k}},\]
where $\alpha_{i_{j}}$ replaces $\partial_{i_j}$.
\end{definition}
The interior product by a covariant tensor of degree $l$, $\alpha^{1}\vee\cdots\vee\alpha^{l}$, is then equal to ${\rm i}(\alpha^1)\circ\cdots\circ {\rm i}(\alpha^l)$.
\begin{definition}
The covariant derivative with respect to $\nabla$ of a symbol $S = \phi| Dx|^{\delta} \otimes \partial_{i_1}\vee\cdots\vee\partial_{i_{k}}$ in the direction of a supervector f\/ield $X$ is def\/ined by
\begin{gather*}
\nabla_{X}(S)
  =   \nabla_{X} (\phi| Dx|^{\delta}) \otimes \partial_{i_1}\vee\cdots\vee\partial_{i_{k}}
\\
 \phantom{\nabla_{X}(S)=}{} +\sum_{j=1}^{k} (-1)^{p(X)(p(\phi)+p(i_{1})+\cdots+p(i_{j-1}))} \phi| Dx|^{\delta} \otimes \partial_{i_1}\vee\cdots\vee\overset{(j)}{\nabla_{X}\partial_{i_j}}\vee\cdots\vee\partial_{i_{k}} .
\end{gather*}
\end{definition}
\begin{definition}
The operator of \emph{divergence} with respect to $\nabla$ is the map
\[
{\rm Div}: \ \mathcal{S}_{\delta}\to\mathcal{S}_{\delta}: \ S\mapsto\sum_{j=1}^{n+m}(-1)^{p(j)}{\rm i}(dx^j)\nabla_{\partial_{j}}S .
\]
\end{definition}
We can easily check that this def\/inition does not depend on the chosen coordinate system.

{\samepage \begin{proposition}
If $j \in {\mathbb{N}}$, $l\in{\mathbb{N}}$ and if $A\in\mathcal{S}_{\delta}^{j}(M)$, then we have
\begin{gather*}
\widetilde{{\rm Div}}\big(A^{h} \vee {\mathcal E}^{l}\big)
  =   ({\rm Div} A)^{h}\vee{\mathcal E}^{l}+2(n-m+1)({\rm i}({\rm r})A)^{h}\vee{\mathcal E}^{l+1}
 -l\gamma_{2j+l}A^{h}\vee{\mathcal E}^{l-1} ,
\end{gather*}
where $\widetilde{{\rm Div}}$ stands for the divergence operator with respect to the lifted supercon\-nec\-tion~$\tilde{\nabla}$ on~$\tilde{M}$ and~${\rm r}$ is the tensor introduced in Definition~{\rm \ref{r}}.
\end{proposition}

\begin{proof}
The proof is exactly the same as in \cite{Bor}.
\end{proof}}

If $\delta$ is not critical, the restriction of $\Psi$ to the divergence-free tensors with respect to~$\tilde{\nabla}$ is thus a~bijection. Indeed, the condition of zero divergence allows to f\/ix the symbols $A_{k-j}$ in (\ref{decomp}). These symbols are given by the following equations (for  $0<l<k$):
\begin{gather*}
A_{k-1} = \frac {1}{\gamma_{2k-1}}{\rm Div}\,A_{k}, \\
A_{k-(l+1)} = \frac {1}{(l+1)(\gamma_{2k-(l+1)})}\left({\rm Div}\, A_{k-l}+2(n-m+1){\rm i}({\rm r})A_{k-(l-1)}\right).
\end{gather*}

Finally, the lift of a symbol $S$, denoted by $\tilde{S}$, is obtained by applying to $S$ the inverse of this bijection. This lift is natural thanks to the naturality of the condition linked to the divergence operator. The lift is also projectively invariant thanks to the fact that this condition depends only on $\tilde{\nabla}$.

\section{Construction of the projectively invariant quantization}

\begin{definition}
If $T$ is a supersymmetric covariant tensor of degree $l$ with values in $\lambda$-densities, $\nabla_{s}T$ is the supersymmetric covariant tensor of degree $l+1$ with values in the $\lambda$-densities def\/ined in the following way:
\begin{gather*}
(\nabla_{s}T)(X_{1},\ldots,X_{l+1}) = \sum_{\sigma\in S_{l+1}}(-1)^{\epsilon_{l+1}+p(T)p(X_{\sigma(1)})}(\nabla_{X_{\sigma(1)}}(T(X_{\sigma(2)},\ldots,X_{\sigma(l+1)}))
\\
\qquad{} -\sum_{j=2}^{l+1}(-1)^{p(X_{\sigma(1)})(p(X_{\sigma(2)})+\cdots+p(X_{\sigma(j-1)}))} T(X_{\sigma(2)},\ldots,\nabla_{X_{\sigma(1)}}X_{\sigma(j)},\ldots,X_{\sigma(l+1)})) ,
\end{gather*}
where $X_{1},\ldots,X_{l+1}$ are supervector f\/ields and where $\epsilon_{l+1}$ is the sign of the permutation $\sigma'$ induced by $\sigma$ on the ordered subset of all odd elements among $X_{1},\ldots,X_{l+1}$.
\end{definition}

\begin{definition}
If $\varphi X_{1}\vee\cdots\vee X_{k}$ is a supersymmetric contravariant tensor of degree $k$ and if $\psi\alpha_{1}\vee\cdots\vee\alpha_{k}$ is a supersymmetric covariant tensor of degree $k$, then we set
\begin{gather*}
\langle\varphi X_{1}\vee\cdots\vee X_{k},\psi\alpha_{1}\vee\cdots\vee\alpha_{k}\rangle\\
\qquad{} =
 \varphi\psi (-1)^{p(\psi)(p(X_{1})+\cdots+p(X_{k}))}{\rm i}(X_{1})\cdots {\rm i}(X_{k})(\alpha_{1}\vee\cdots\vee\alpha_{k}),
 \end{gather*}
where the interior product ${\rm i}$ is def\/ined in the same way as in Def\/inition~\ref{interior}. One extends this operation by bilinearity to arbitrary supersymmetric tensors of degree $k$.
\end{definition}

\subsection{The main result}

In this section, we give an explicit formula for the natural and projectively invariant quantization.
\begin{theorem}\label{princ} If $n-m \neq \pm 1$  and $\delta$ is not a critical value, then the collection of maps $Q_M: \mathcal{C}\times \mathcal{S}_{\delta} \to \mathcal{D}_{\lambda,\mu}$ given by
\begin{gather}\label{formula}
{\left( Q_M(\nabla,S)(f) \right)}^{\sim} = \langle\tilde{S} , \tilde{\nabla}_s^{k} \tilde{f}\rangle,
\end{gather}
defines a projectively invariant natural quantization for supermanifolds of dimension~$(n|m)$.
\end{theorem}

\begin{proof}
The proof goes exactly in the same way as in~\cite{Bor}. First, formula (\ref{formula}) is well-def\/ined: indeed, the right-hand side is $\mu$-equivariant because of the invariance of $\tilde{\nabla}$, of the $\delta$-equivariance of $\tilde{S}$ and the $\lambda$-equivariance of $\tilde{f}$.

Next, (\ref{formula}) def\/ines obviously a natural and projectively invariant quantization: this fact is mainly due to the naturality and the projective invariance of the lift $\tilde{S}$. The quantization preserves the principal symbol for the same reasons as in~\cite{Bor}.
\end{proof}

\begin{remark}\label{rem:cas superplat n-m=-1}
When $n-m \neq \pm 1$, $M = {\mathbb{R}}^{n| m}$ and $\nabla = \nabla_0$, formula (\ref{formula}) recovers the unique $\mathfrak{pgl}(n+1|m)$-equivariant quantization found in \cite{MR5}. It is interesting to notice the problem there was solved without any hypothesis on the superdimension.
\end{remark}

\begin{remark}
When $n-m \neq -1, 1, -2, -4$ and $\delta = 0$, formula (\ref{formula}) coincides with the canonical dif\/ferential operator associated by J.~George \cite[Theorem 3.6]{George} with a symbol of degree two and a projective class of superconnections. In particular, when $\lambda=\mu=0$, our formula recovers the projective Laplacian of \cite[Theorem 3.2]{George}. Notice that J.~George does not have any formula for $n-m \in \{ -2, -4 \}$ while we do in general (those additional conditions appear only for particular values of the shift $\delta$). Actually, the methods used are dif\/ferent. The conditions $n-m \neq -1,1,-2$ in his work come from the construction, given a projective class on $M$, of the quantities $\tilde{\Pi}_{\mathfrak{b}\mathfrak{c}}^{\mathfrak{a}}$ associated with the corresponding Thomas connection while the condition $n-m \neq -4$ is due to the use of the projective Laplacian on $\tilde{M}$ (see \cite[Theorem 3.2]{George} for details). This being said, the relation between projectively invariant quantization and J.~George's work is still unclear: one can wonder how the procedures are linked in the general case~$\delta \neq 0$.
\end{remark}

\section[The case $n-m = 1$]{The case $\boldsymbol{n-m = 1}$}

As it was already noticed in Section~\ref{subsection:Thomas connection}, the hypothesis $n-m \neq 1$ is the analogue in the super context of the fact that M.~Bordemann's method~\cite{Bor} does not apply for $1$-dimensional smooth manifolds.

Actually, the problem of natural and projectively invariant quantization on $1$-dimensional smooth manifolds turns out to be very peculiar. In this case, it is easily shown that the dif\/ference between any two torsion-free linear connections can be expressed as $\alpha \vee {\mathrm{id}}$ for some $1$-form $\alpha$. Consequently, all torsion-free linear connections are projectively equivalent, and the quest for a natural projectively invariant quantization amounts to the quest for a natural bijection from symbols to dif\/ferential operators. As it is well-known (it is for instance a consequence of~\cite[Theorem 3]{gargoubi}), such a natural bijection does not exist. Notice that for symbols of order two, the theory of natural operators~\cite{KMS} imposes for a natural projectively invariant quantization to be of the form
\begin{gather}\label{eq:natural quantization}
Q(\nabla,S)(f) = \langle S,\nabla^2 f \rangle + a   \langle {\rm Div}\, S,\nabla f \rangle + b   \langle {\rm Div^2}S,f \rangle + c   \langle{\rm i}({\rm Ric}) S,f \rangle,
\end{gather}
where $a,b,c \in {\mathbb R}$. The condition of projective invariance yields a system of equations for $a$, $b$, $c$ which admits no solution in dimension $n=1$ (cf.~\cite{bouarroudj}).

If we make the assumption that a natural projectively invariant quantization must write under the form~(\ref{eq:natural quantization}), with all objects being replaced by their super analogues, then the system of equations provided by the condition of projective invariance has no solutions when $n-m = 1$. Therefore, unless there are more natural operators for supermanifolds than the superizations of the classical ones, a natural projectively invariant quantization does not exist in this case.

\section[The case $n-m = -1$]{The case $\boldsymbol{n-m = -1}$}

In \cite{MR5}, P.~Mathonet and F.~Radoux were able to build a $\mathfrak{pgl}(n+1|m)$-equivariant quantization without any hypothesis on the superdimension. Nevertheless, the case $n-m = -1$ required an ad-hoc construction because of the peculiarities of the Lie superalgebra $\mathfrak{pgl}(n+1|n+1)$.

In our case, the problem lies in the very def\/inition of the quantities $\Pi^i_{jk}$ used in the construction of the connection $\tilde{\nabla}$ on $\tilde{M}$ associated with a projective class of connections on~$M$. The manifold $\tilde{M}$ is thus unhelpful here.

This being said, it can be checked by hand that the formula
\[
Q(\nabla,S)(f) =  \langle S,\nabla f \rangle + t    \langle {\rm Div}\,S, f \rangle
\]
def\/ines a $1$-parameter family of natural projectively invariant quantization for symbols of order one. This result agrees with the phenomenon observed in~\cite{MR5}. Also, for symbols of order two, the formula
\[
Q(\nabla,S)(f) =  \langle S,\nabla^2 f \rangle + \langle {\rm Div}\,S, \nabla f \rangle
\]
turns out to be projectively invariant. We conjecture that similar formulas can be obtained for higher order symbols and that a natural projectively invariant quantization exists when $n-m = -1$.

\subsection*{Acknowledgements}

It is a pleasure to thank P. Mathonet for fruitful discussions. We also thank the referees for suggestions leading to great improvements of the original paper. Finally, F.~Radoux thanks the Belgian FNRS for his research fellowship.

\pdfbookmark[1]{References}{ref}
\LastPageEnding

\end{document}